\documentclass[review]{elsarticle}
\usepackage[symbol]{footmisc}
 \usepackage[intlimits]{amsmath}
\usepackage{lineno,hyperref}
\usepackage{float}
\usepackage{pgfplots}
\pgfplotsset{compat=newest}
\usepackage{tikz}
\usetikzlibrary{patterns}
\usepackage{rotating}
\usepackage{graphicx} 
\usepackage{graphics}
\usepackage{standalone} 
\usepackage{cleveref}
\input{IAM_Lager.sty}
\usetikzlibrary{snakes}
\usetikzlibrary{calc}
\usetikzlibrary{arrows,matrix,positioning}
\usetikzlibrary{fit}
\usetikzlibrary{tikzmark}
\usepackage{eso-pic}
\usepackage{xcolor}
\usepackage{array}
\usepackage{blkarray}
\usepackage{setspace} 
\usepackage{booktabs}
\usepackage{array}
\usepackage{cellspace}
\usepackage{makecell}
\usepackage{caption}
\usepackage{subcaption}
\usepackage{multirow}
\usepackage{longtable}
\usepackage{comment}
\usepackage{amssymb}
\usepackage{tablefootnote}
\usepackage{multirow}

\modulolinenumbers[5]

\newcommand{\pur}[1]{\textcolor{gray}{#1}}
\newcommand{\pos}{{\xi}}
\newcommand{\setNodes}{{\mathcal{N}}}
\newcommand{\setBeams}{{\mathcal{T}}}
\newcommand{\setBearings}{{\mathcal{Q}}}
\newcommand{\setSprings}{{\mathcal{S}}}
\newcommand{\setLoads}{{\mathcal{F}}}
\renewcommand{\vec}[1]{\mathbf{#1}}

\setcellgapes{4pt}

\usetikzlibrary{calc}
\def\centerarc[#1](#2)(#3:#4:#5)
{ \draw[#1] ($(#2)+({#5*cos(#3)},{#5*sin(#3)})$) arc (#3:#4:#5); }

\crefformat{equation}{Eq.~(#2\textup{#1}#3)}
\tikzset{
	hatch distance/.store in=\hatchdistance,
	hatch distance=7pt,
	hatch thickness/.store in=\hatchthickness,
	hatch thickness=0.5pt
}

\makeatletter
\pgfdeclarepatternformonly[\hatchdistance,\hatchthickness]{flexible hatch}
{\pgfqpoint{0pt}{0pt}}
{\pgfqpoint{\hatchdistance}{\hatchdistance}}
{\pgfpoint{\hatchdistance-1pt}{\hatchdistance-1pt}}%
{
	\pgfsetcolor{\tikz@pattern@color}
	\pgfsetlinewidth{\hatchthickness}
	\pgfpathmoveto{\pgfqpoint{0pt}{0pt}}
	\pgfpathlineto{\pgfqpoint{\hatchdistance}{\hatchdistance}}
	\pgfusepath{stroke}
}
\makeatother

\journal{Journal of Sound and Vibration}

\begin{document}

\begin{frontmatter}
\title{The Numerical Assembly Technique for arbitrary planar systems based on an alternative homogeneous solution}
\author{Thomas Kramer\textsuperscript{\textit{a,}}\footnote[1]{kramer@tugraz.at} , Michael Helmut Gfrerer\textsuperscript{\textit{a}}  
	\vspace*{5mm}}
\address{\textsuperscript{\textit{a}}Institute of Applied Mechanics 
	\\ \vspace*{2mm} University of Technology Graz \\ \vspace*{2mm}
	Technikerstraße 4, 8010 Graz, Austria}

\begin{abstract}	
The Numerical Assembly Technique is extended to investigate arbitrary planar frame structures with the focus on the computation of natural frequencies. This allows us to obtain highly accurate results without resorting to spatial discretization. 
To this end, we systematically introduce a frame structure as a set of nodes, beams, bearings, springs, and external loads and formulate the corresponding boundary and interface conditions. As the underlying homogeneous solution of the governing equations, we use a novel approach recently presented in the literature. This greatly improves the numerical stability and allows the stable computation of very high natural frequencies accurately. We show this numerically at two frame structures by investigation of the condition number of the system matrix and also by the use of variable precision arithmetic. 	
\end{abstract}

\end{frontmatter}
\section{Introduction}
A significantly high number of mechanical and structural engineering problems can be modeled by an assembly of beams with attached supports, hinges, and springs such as building structures but also robotic arms or any mechanical mechanism. Especially for structural features where the dynamic behavior is of major interest, the natural frequencies play an important role. Thus, having an efficient tool to determine those natural frequencies, combined with the corresponding mode shapes, helps engineers perform high-quality work.

For simple problems analytical formulas can be found in various textbooks \cite{dinkler2016einfuhrung}-\nocite{rao2019vibration,gasch1989strukturdynamik}\cite{schiehlen2017technische}. There, functions to compute natural frequencies and mode shapes are derived for e.g. cantilever beams or single-span beams with various supports. Also for more complicated systems a huge amount of analytical and semi-analytical methods can be found in the literature. A detailed literature review can be found in \cite{klanner2020steady}.
 
One widely used method is the Dynamic Stiffness Method (DSM) which uses shape functions tailored to the governing equations and the assembly technique from the Finite Element Method \cite{mario2019structural}. Due to its general structure, it can be applied to arbitrary frame structures. Recently, it has been extended to functionally graded beams \cite{banerjee2018free} and axially loaded multi-cracked frames \cite{caddemi2009exact} among other extensions such as beams carrying spring-mass systems \cite{banerjee2012free} and axially loaded composite Timoshenko beams \cite{banerjee1998free}.

To analyze multi-stepped beams with various concentrated masses the Transfer Matrix Method (TMM) is often used. With this approach, a matrix transfers the beam state such as the displacement, the rotation, the internal forces, and the bending moment from one end of a uniform beam segment to the other. Anyway, this method's application is limited to multi-span \cite{bapat1987natural} or multi-stepped beams \cite{el2017normalized} and can, therefore, not be applied on plane frame structures. Also, this method has been extended to crack identification of stepped beams with multiple edge cracks and different boundary conditions \cite{attar2012transfer}.

Another noteworthy method is the Green Function Method (GFM), which uses the analytical responses of a uniform beam segment to point loads (Green functions) to include concentrated elements in the beam model. It is left to solve for the displacement and the rotations of the beam at each element. Hence, the size of the system matrix depends clearly on the number of concentrated elements. Lately, this method has been used to analyze elastically supported Euler-Bernoulli beams \cite{ronvcevic2019closed}, cracked beams with elastic boundary conditions \cite{ghannadiasl2018analytical} and axially functionally graded and non-uniformed beams \cite{han2019new}.

The Generalized Function Approach (GFA) follows a similar idea as the GFM, where the discontinuous response of the beam due to concentrated elements and the changes in material parameters is handled with generalized functions such as Dirac delta functions or Heaviside step functions. Contrarily to the GFM, the size of the system matrix is independent of the number of uniform beam sections. A plane vibration problem can, therefore, be assembled by a $2\times 2$ system matrix.
The GFA has been extensively studied. Most recently on beams with various in-span supports and masses \cite{burlon2018exact}, on discontinuous layered elastically bonded beams \cite{di2018flexural}, on plane structures with mass-spring subsystems and rotational joints \cite{failla2019exact} and on beams including tuned mass
dampers with spring inertia effects \cite{failla2019random}.  

The Numerical Assembly Technique (NAT) considers the analytical solution of the homogeneous governing equation. Each analytical solution with the enforced boundary and intermediate conditions holds within a uniform beam segment and occupies one row in the system of equations. Hence, the implementation is straight forward and the size of the system of equations rises with the complexity of the structure.

Numerous research papers concerning the NAT have been published. Whereas two authors have caught our attention specifically throughout the literature studies as they are certainly experts in this field of research. Their most recent publications regarding the NAT are the computation of steady-state harmonic vibrations of viscoelastic Timoshenko Beams with fractional derivative damping models \cite{applmech2040046}, the application to multi-stepped Euler-Bernoulli beams under arbitrarily distributed loads carrying any number of concentrated elements \cite{klanner2020steady}, the free vibration analysis of continuous bridges under vehicles as ambient excitation forces \cite{tan2017free2} and the application to cracked Timoshenko beams carrying spring-mass systems \cite{tan2017free}.
However, to the best knowledge of the authors, the NAT is only applied to multi-span and/or multi-stepped beams.

In the present paper, the NAT is extended to plane frame structures based on the Euler-Bernoulli beam theory with various combinative supports, hinges, and springs. Therefore, we systematically introduce a frame structure as a set of nodes, beams, bearings, springs, and external loads and formulate the corresponding boundary and interface conditions. Writing the system of equations in matrix form allows us to solve for the natural frequencies of the plane structure. This has to be done numerically.
More than that, the general function to solve the governing equation for transversal vibrations is improved by using the novel approach introduced in \cite{klanner2020steady}. This increases the accuracy of the results, especially when computing natural frequencies in high-frequency ranges. The proposed NAT is implemented in a \textsc{Matlab}\textcopyright-based code and its numerical performance is investigated at two exemplary frame structures.

\section{Mechanical Models}\label{sec:mechanicalModels}
\subsection{Governing equations for a single beam}
When considering vibrations in the two-dimensional space, the longitudinal and transversal deformation of an element are implicated. Let $\ell$ be the length of the considered beam with constant geometrical and material parameters along its local beam axis, where $\rho$ is the density, $E$ the Youngs Modulus, $A$ the cross-section area, and $I$ the area moment of inertia. With the position $\pos \in[0,\ell]$ and the time $t$, the differential equations for longitudinal vibrations are given by
\begin{subequations}
	\begin{align}
	\rho  A \ddot{u}(\pos,t) &= N'(\pos,t) + n(\pos,t), \label{eq:equCondLong} \\
		N(\pos,t) &= EA  u'(\pos,t).
	\end{align}
\end{subequations}
With $u(\pos,t)$ as the axial deformation, $N(\pos,t)$ defining the internal axial force and $n(\pos,t)$ the external axial loading force. Further, $()'$ denotes the derivative with respect to $\xi$ and $(\dot)$ the derivative with respect to time.
The differential equations for transversal vibrations governed by the Euler-Bernoulli beam theory read as
\begin{subequations}
	\begin{align}
	\rho A \ddot{w}(\pos,t) &= V'(\pos,t)+q(\pos,t), \\
	V(\pos,t) &= M'(\pos,t),\\
	M(\pos,t) &= EI\psi'(\pos,t),\\
	w'(\pos,t) &=-\psi(\pos,t).
	\end{align}
\end{subequations}
Where the transversal deflection $w(\pos,t)$ is described through the internal shear force $V(\pos,t)$, the internal bending moment $M(\pos,t)$, the rotation angle $\psi(\pos,t)$ and the external transversal loading force $q(\pos,t)$.
When referring to the natural frequencies of a structure, the considered equation simplifies to a homogeneous differential equation, since it holds for the load forces $n=0$ and $q=0$. 
Assuming free harmonic vibrational behavior with the natural frequency $\omega$ the governing equations are
	\begin{equation}
	\frac{\partial^2 U(\pos)}{\partial \pos^2} + c^2 \ U(\pos) = 0,	
\label{eq:diffEqLong}
\end{equation}
\begin{equation}
	\frac{\partial^4 W(\pos)}{\partial \pos^4} - \kappa^4 \ W(\pos) = 0,
	\label{eq:diffEqTravers}
\end{equation}
with
\begin{align*}
c^2 = \frac{\rho}{ E}\omega^2, \qquad
\kappa^4 =  \frac{\rho A}{ EI}\omega^2.
\end{align*}
\subsection{Frame structure}\label{sec:frame}
To introduce a frame structure formally, we extend the notation introduced in \cite{polz2019wave}. Let $\setNodes$ be the set of nodes, $\setBeams$ the set of beams, $\setBearings$ the set of bearings, $\setSprings$ the set of springs and $\setLoads$ the set of external loads. 
For a frame structure $(\setNodes,\setBeams,\setBearings,\setSprings,\setLoads)$ these sets have the following properties:
\begin{enumerate}
	\item Each beam $m \in \setBeams$ connects two nodes $k_0,k_\ell \in \setNodes$ and holds constant parameters $\rho_m,\, E_m,\, I_m,\, A_m$. The beam $m$ can either be connected rigidly or hinged at $k_0$ and $k_\ell$.
	\item Each node $k \in \setNodes$ has a set of attached beams $\beta_k \subseteq \setBeams$, a set of attached springs $\chi_k \subseteq  \setSprings$, external loads $\vec{F}_k^{ext}\in\setLoads$ and external moments $M_k^{ext}\in\setLoads$. Further, a bearing $\gamma_k \subseteq \setBearings$ can be attached at any node $k$.
\end{enumerate}
Here, we assumed that the external forces and moments only act on the nodes. In general, also distributed loads acting on the beams are possible. However, due to the focus on the computation of natural frequencies, the influence of external loads is not further considered in \cref{eq:diffEqLong} and \cref{eq:diffEqTravers}.
For a beam $m \in \setBeams$ we define the tangent vector $\vec{t}_m =k_{\ell,m}- k_{0,m}$ and the normal vector $\vec{n}_m\perp\vec{t}_m$. The resulting internal force vector of beam $m\in\setBeams$ is defined as
\begin{equation}
\vec{S}_m(\pos)=N_m(\xi)\cdot\vec{t}_m+V_m(\xi)\cdot\vec{n}_m.
\end{equation}
\subsection{Boundary and coupling conditions}
With the definitions in \Cref{sec:frame} the boundary and interface conditions can now be formulated. First, the kinetics are defined by formulating the equilibrium of forces at each node $k$,
\begin{equation}
	\vec{F}_k = \sum_{m \in \beta_k}^{} \vec{S}_{m|k}+\vec{Q}_k^L+\vec{F}_k^{ext}+\vec{F}_k^L=0,
	\label{eq:equilibrium_F}
\end{equation}
where,
	$\sum \vec{S}_{m|k}$ is the sum of each beam's internal forces,
	$\vec{Q}_k^L$ the reaction forces from bearings,
	$\vec{F}_k^{ext}$ the external loading forces and
	$\vec{F}_k^L$ the constrained forces due to springs.
	The equilibrium of moments at each node $k$ reads
\begin{equation}
	{M}_k = \sum_{m \in \beta_k}^{} {M}_{m|k}+{Q}_k^R+{M}_k^{ext}+{M}_k^R=0,
	\label{eq:equilibrium_M}	
\end{equation}
where,
	$\sum {M}_{m|k}$ is the sum of each beam's internal bending moment,
	${Q}_k^R$ the constrained rotational moments from bearings,
	${M}_k^{ext}$ the external loading moments and
	${M}_k^R$ the rotational moments due to springs.
	The terms in \cref{eq:equilibrium_F} and \cref{eq:equilibrium_M} that include springs are defined through Hooke's law,
	\begin{equation}
		\vec{F}_k^L=c_L(\vec{U}_k\cdot\pmb\zeta)\cdot\pmb\zeta
	\end{equation}
and
	\begin{equation}
		{M}_k^R=c_\psi\cdot\psi_k	,
	\end{equation}
	where $c_L$ and $c_\psi$ define the spring constants of the longitudinal and the rotational spring, respectively and $\pmb\zeta$ is the direction vector along the axis of the longitudinal spring. 
	Ultimately $\vec{U}_k$ and $\psi_k$ describe the displacement and the rotation of node $k$, respectively.
Next we formulate the coupling conditions for a node $k \in \setNodes : |\beta_k|>1$. Therefore, let $m_f \in \beta_k$ be a beam that is fixed in the set of beams connected at $k$. The compatibility condition of the resulting deformation $\vec{U}$ of the beam $m$ evaluated at the node $k$ reads 

\begin{equation}
	\vec{U}_{m_f|k}=\vec{U}_{m|k} \qquad \forall \ m \in \beta_{k} \backslash m_f.
	\label{eq:coupling1}
\end{equation}

To differentiate between rigid and hinged connections the set of attached beams $\beta_k$ at a node $k \in \setNodes : |\beta_k|>1$ is divided into two subsets, so that $\beta_k=\beta_k^r \cup \beta_k^h$. Thus, $\beta_k^r$ includes the beams that are rigidly connected with one another and $\beta_k^h$ the beams that are attached at $k$ with an hinge. 
To formulate the kinematic relations concerning rotations, let $m_r\in \beta_k^r$ be a rigidly connected beam at node $k$ that is fixed in the set of rigidly connected beams at this node. For all $m \in \beta_k^r \backslash m_r$ the kinematic conditions for the rotation $\psi$ of the beam $m$ evaluated at the node $k$ are
\begin{equation}
	\psi_{m_r|k}=\psi_{m|k} \qquad \forall \ m \in \beta_{k}^r \backslash m_r.
	\label{eq:coupling2}
\end{equation}
Furthermore, let $m_h\in \beta_k^h$ be a beam that is connected with a hinge at the node $k$ and fixed in the set of beams with a hinged connection at this node. The kinetic conditions are

\begin{equation}
	M_{m|k}=0 \qquad \forall \ m \in \beta_{k}^h\backslash m_h.
	\label{eq:coupling3}	
\end{equation}
It is noteworthy that through the equilibrium equation \cref{eq:equilibrium_M} the condition $M_{m_h|k}=0$ holds true.
The boundary conditions due to applied bearings are defined in the local coordinate system for bearings $[\zeta,\eta]$ with the corresponding unit vectors $\vec{e_\zeta}\perp\vec{e_\eta} $, see \Cref{tab:supports} and read for the pinned bearing
\begin{align}
	\vec{U}_{m|k}=0 \qquad  \ m \in \beta_{k},
	\label{eq:support1}	
\end{align}
the roller
\begin{align}
	 \vec{U}_{m|k}\cdot \vec{e_\eta}=0  \qquad  \ m \in \beta_{k},
	\label{eq:support2}		 
\end{align}
the clamped bearing
\begin{subequations}
\begin{align}
	\vec{U}_{m|k}=0 \qquad  \ m \in \beta_{k},\label{eq:support3a}\\
	 \psi_{m|k}=0 \qquad  \ m \in \beta_{k}^r,\label{eq:support3b}			 
\end{align}
\end{subequations}
and the parallel guide
\begin{subequations}
\begin{align}
	\vec{U}_{m|k}\cdot\vec{e_\eta}=0 \qquad  \ m \in \beta_{k},\label{eq:support4a}\\
	\psi_{m|k}=0 \qquad  \ m \in \beta_{k}^r,
	\label{eq:support4b}		
\end{align}
\end{subequations}
respectively.
\section{Local Solutions of the governing equations}
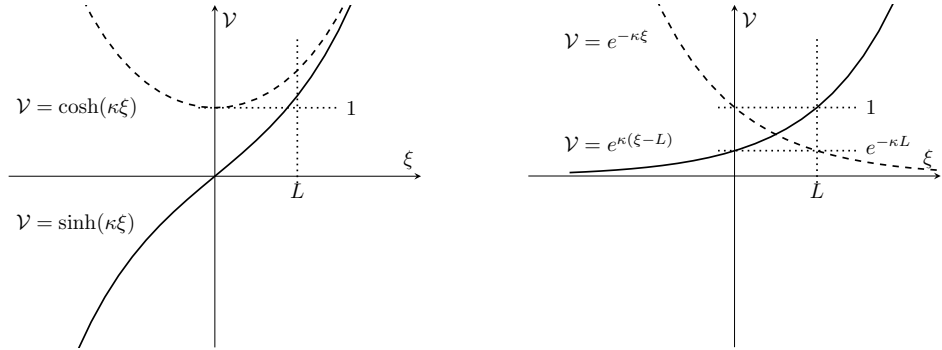
\begin{figure}[H]
	\begin{minipage}{0.43\textwidth}
		\scalebox{.8}{\begin{tikzpicture}
\begin{axis}[
    xmin=-2.5, xmax=2.5,
    ymin=-2.5, ymax=2.5,
    axis lines=center,
    axis on top=true,
    domain=-2:2,
    ylabel=$\mathcal{V}$,
    xlabel=$\xi$,
    xmajorticks=false,
    ticks=none
    ]

    \addplot [mark=none, thick] {sinh(\x)};
   	\node [right] at (axis cs: -2.5,-0.7) {$\mathcal{V} = \sinh (\kappa \xi)$};
   	
   	\addplot [mark=none,dashed, thick] {cosh(\x)};
   	\node [right] at (axis cs: -2.5,1) {$\mathcal{V} = \cosh (\kappa \xi)$};
    \draw [dotted, thick] (axis cs:1,-0.1)-- (axis cs:1,2);
    \node [below] at (axis cs: 1,0) {$L$};
    
    \draw[dotted,thick](axis cs:-.2,1)-- (axis cs:1.5,1);
    \node [right] at (axis cs:1.5,1){1};
\end{axis}
\end{tikzpicture}}
	\end{minipage}
	\hfill
	\begin{minipage}{0.43\textwidth}
		\scalebox{.8}{\begin{tikzpicture}
\begin{axis}[
    xmin=-2.5, xmax=2.5,
    ymin=-2.5, ymax=2.5,
    axis lines=center,
    axis on top=true,
    domain=-2:2.5,
    ylabel=$\mathcal{V}$,
    xlabel=$\xi$,
    xmajorticks=false,
    ticks=none
    ]

    \addplot [mark=none, thick] {exp(\x-1)};
   	\node [right] at (axis cs: -2.155,0.5) {$\mathcal{V} =  e^{\kappa(\xi-L)}$};
   	
   	\addplot [mark=none,dashed, thick] {exp(-\x)};
   	\node [right] at (axis cs: -2.155,2) {$\mathcal{V} = e^{-\kappa \xi}$};
    \draw [dotted, thick] (axis cs:1,-0.1)-- (axis cs:1,2);
    \node [below] at (axis cs: 1,0) {$L$};
    
    \draw[dotted,thick](axis cs:-.2,1)-- (axis cs:1.5,1);
    \node [right] at (axis cs:1.5,1){1};
    \draw[dotted,thick](axis cs:-.2,.37)-- (axis cs:1.5,.37);    
    \node [right] at (axis cs:1.5,.45){$e^{-\kappa L}$};
\end{axis}
\end{tikzpicture}}
	\end{minipage}
	\caption{The difference of ansatz functions for the general solution. The hyperbolic- (left) and exponential approach (right).}
	\label{fig:hyperExp}
\end{figure}

\begin{figure}[H]
	\begin{minipage}{.39\textwidth}
		\begin{tikzpicture}
			\begin{scope}[scale=3]
				\def\lBeam{1}
				\def\offset{\lBeam/25}
				\def\newLine{.2}
				\coordinate (A) at (0,0);
				\coordinate (B) at (\lBeam,0);
				\draw[thick](A)--(B);
				\draw[dashed]($(A)+(\offset,-\offset)$)--($(B)+(-\offset,-\offset)$);
				\Festlager{A}[0][.3];
				\Festlager{B}[0][.3];
				\Bemassung{A}{B}{L}[.3];
				\coordinate (Text) at ($(A)+(.5*\lBeam,-.4*\lBeam)$);				
				\node[align=center] at ($(Text)$) {\small$L=1\ [m]$};
				\node[align=center] at ($(Text)+(0,-\newLine)$) {\small$A= 7.56 \cdot 10^{-4} \ [m^2]$};
				\node[align=center] at ($(Text)+(0,-2*\newLine)$) {\small$ I = 3.5 \cdot 10^{-10} \ [m^4] $};
				\node[align=center] at ($(Text)+(0,-3*\newLine)$) {\small$ \rho =  7.85 \ [t/m^3]$};
				\node[align=center] at ($(Text)+(0,-4*\newLine)$) {\small$ E = 2 \cdot 10^{8} \ [kN/m^2]$};	
			\end{scope}
		\end{tikzpicture}
	\end{minipage}
	\begin{minipage}{.59\textwidth}
	 	\input{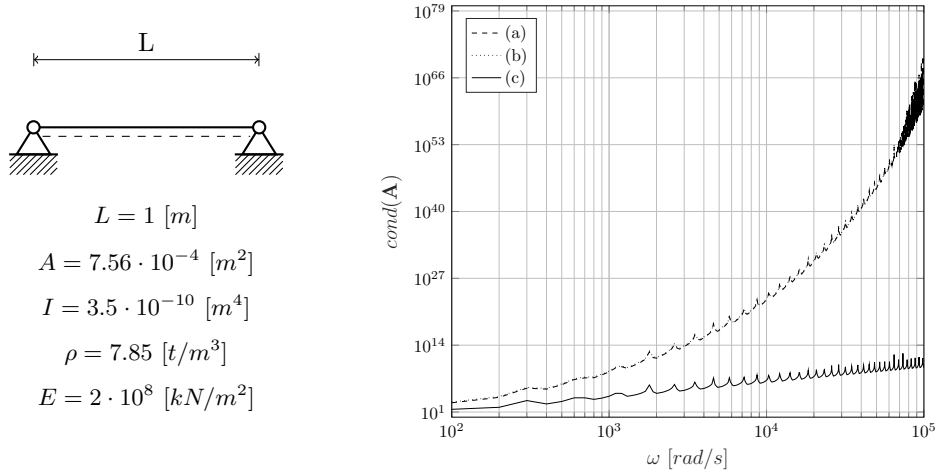}	
	\end{minipage}
\caption{The condition number for a single span beam's system matrix.}
\label{fig:condNum}
\end{figure}

The local solution of the differential equation \cref{eq:diffEqLong} in the frequency domain reads as 
\begin{equation}
U_m{(\pos)} = C_{m,1} \cos\left(c_m\pos\right) + C_{m,2}\sin\left(c_m\pos\right).
\label{for:general_sol_SR}
\end{equation}
For the transversal vibrations, in this work, we refer to the \textit{alternative ansatz} which differs from the conventionally used solutions. Therefore, its derivation is stated.
When implementing the ansatz $W_{(\pos)} = e^{\lambda \pos}$, the equation \cref{eq:diffEqTravers} can be expressed with its characteristic polynomial
\begin{equation}
\lambda^4-\kappa^4=0.
\label{eq:diffEq2}
\end{equation}
The four roots of the eigenvalue problem \cref{eq:diffEq2} are $\lambda_1=+\kappa$, $\lambda_2=+i\kappa$, $\lambda_3=-\kappa$, and $\lambda_4=-i\kappa$. Thus, a general solution that corresponds the real and complex roots of \cref{eq:diffEq2} needs to be found. Many publications refer to the commonly used ansatz
\begin{equation}
	W_m{(\pos)} = C_{m,3}\cos(\kappa \pos)+C_{m,4}\sin(\kappa \pos)+C_{m,5} \cosh(\kappa \pos)+ C_{m,6}\sinh(\kappa \pos).
	\label{for:eigenfunction}
\end{equation}
In this publication, a different ansatz function for the general solution of the differential equation for transversal vibration is presented. When rewriting the ansatz $W{(\pos)} = e^{\lambda \pos}$ with the obtained roots $\lambda_{1-4}$, $W(\pos)$ can be expressed as
\begin{equation}
	W_m{(\pos)} = C_{m,1} e^{+i\kappa \pos} + C_{m,2} e^{-i\kappa \pos} + C_{m,3} e^{+\kappa \pos} + C_{m,4} e^{-\kappa \pos}.
	\label{eq:alternAnsatz01}
\end{equation}
Since the condition
\begin{equation}
	\cos(\pos)+sin(\pos)=\left(\frac{1}{2} + \frac{i}{2}\right) e^{-i \pos} + \left( \frac{1}{2} - \frac{i}{2} \right) e^{i \pos}
	\label{eq:condition01}
\end{equation}
holds, the first two terms in \cref{eq:alternAnsatz01} can be simplified as trigonometrical functions, whereas the last two terms in \cref{eq:alternAnsatz01} are kept as exponential functions so that the general solution is 
\begin{equation}
	W_m{(\pos)} = C_{m,1}\cos(\kappa \pos)+C_{m,2}\sin(\kappa \pos)+C_{m,3} e^{\kappa \pos}+ C_{m,4} e^{-\kappa \pos}.
	\label{for:eigenfunction2.1}
\end{equation}
After dividing the third term in \cref{for:eigenfunction2.1} by $e^{\kappa L}$, so that
\begin{equation}
	e^{\kappa \pos} \cdot e^{-\kappa L} = e^{\kappa(\pos-L)},
	\label{eq:condition02}
\end{equation}
the final general solution is
\begin{equation}
	W_m{(\pos)} = C_{m,1}\cos(\kappa \pos)+C_{m,2}\sin(\kappa \pos)+C_{m,3} e^{\kappa(\pos-L)}+ C_{m,4} e^{-\kappa \pos}.
	\label{for:eigenfunction2}
\end{equation}
In the alternative ansatz \cref{for:eigenfunction2}, hyperbolic functions are avoided in the general solution. The reason for that is the asymptotic form of hyperbolic functions. As the value of $\pos$ or $\omega$ increases, the hyperbolic function terms get extremely high, and therefore, handling with those high numbers becomes rather cumbersome and the computational effort increases drastically. Contrarily, the amplitudes of the exponential terms in the general solution \cref{for:eigenfunction2} stay within the fixed interval $\mathcal{V} \in [e^{-\kappa L},1]$, see \Cref{fig:hyperExp}. Hence, determining the natural frequency using this ansatz is much simpler and faster.
To illustrate this, we investigate the condition number for the bending terms of the system matrix $\vec{A}$ of a one-span beam, see \Cref{fig:condNum}.
The condition number gives information about the sensitivity of a linear system to small changes in a matrix.
When using the Numerical Assembly Technique to determine the natural frequencies of a structure, a system matrix is generated. Depending on the number of beams, the size of this matrix can get rather big. Due to rounding errors, small changes in the system matrix could lead to instability of the linear system when choosing the common hyperbolic functions in the general solution. The proposed exponential and trigonometrical approach ensures stable computations even at high frequencies. The graphs in \Cref{fig:condNum} refer to the three different general solutions \cref{for:eigenfunction}$\to$(a), \cref{for:eigenfunction2.1}$\to$(b) and \cref{for:eigenfunction2}$\to$(c).
It clearly shows in \Cref{fig:condNum} that the linear system based on the commonly used general solution (a), as well as the second approach (b), becomes quickly unstable for higher frequencies. The graphs of the solutions (a) and (b) coincide with this structure and rise exponentially in this double-logarithmic plot. But the instability of the linear system using the modified general solution (c) rises linearly in a relatively moderate manner.
Thus, the evaluation of natural frequencies should therefore be improved by this change and has an effect when the natural frequencies need to be accurately determined, particularly for higher frequency ranges.
\section{The Numerical Assembly Technique}
In the Numerical Assembly Technique, the homogeneous solutions \cref{for:eigenfunction}
, \cref{for:eigenfunction2.1}
 and \cref{for:eigenfunction2} 
 as well as \cref{for:general_sol_SR}, are used to fulfill the boundary and interface conditions, stated in \Cref{sec:mechanicalModels}, and therefore, the harmonic governing equations \cref{eq:diffEqLong} and \cref{eq:diffEqTravers} are exactly satisfied. The system of equations is
\begin{equation}
	\mathbf A(\omega) \mathbf c = \mathbf 0
	\label{eq:sysOfEq}	
\end{equation}
where $\mathbf c$ is a $n\times1$ vector $[\vec{c_m},\vec{Q}]^T$ containing the integration constants $\vec{c_m}$ and the bearing forces and moments $\vec{Q}$ and $\vec{A}$ forms the $n\times n$ asymmetric system matrix and only depends on $\omega$. The size of the system of equations is defined as $n=|\setBeams|\cdot6+|\setBearings|$. Each row of $\vec{A}$ is obtained by inserting the homogeneous solutions into the boundary and interface conditions \cref{eq:equilibrium_F}, \cref{eq:equilibrium_M}, \cref{eq:coupling1} to \cref{eq:support2}, \cref{eq:support3a} ,\cref{eq:support3b} ,\cref{eq:support4a} and \cref{eq:support4b}, respectively.  
Non trivial solutions for $ \omega $ fulfill the system of equations \cref{eq:sysOfEq} when its determinant is zero,
\begin{equation}
	f(\omega)=\det(\mathbf{A}) = 0.
	\label{eq:detK}
\end{equation}
Here, $ f(\omega) $ is a continuous function with an infinite number of zero crossings and is often referred to as the eigen-, or frequency function. Each $ \omega $ that fulfills the so-called characteristic equation \cref{eq:detK} is an eigenvalue of the system of equations and, therefore, a natural frequency of the structure. In the present work, the zeros of $f(\omega)$ are determined using the Regular-Falsi method. This method gives accurate results over a rather short computation time and, meets therefore the expectations of this application.
Implementing the calculated $ \omega $-values into our system of equations makes it possible to solve for the integration constants $\vec{c_m}$ and the bearing forces and moments $\vec{Q}$ representing the mode shapes to each natural frequency.
\section{Numerical Results}\label{sec:caseStudies}
\subsection{Numerical Precision}
\begin{figure}[H]
	\centering
	\begin{subfigure}[b]{0.9\textwidth}
		\begin{minipage}{.4\textwidth}
			\centering
			\includegraphics[scale=0.3]{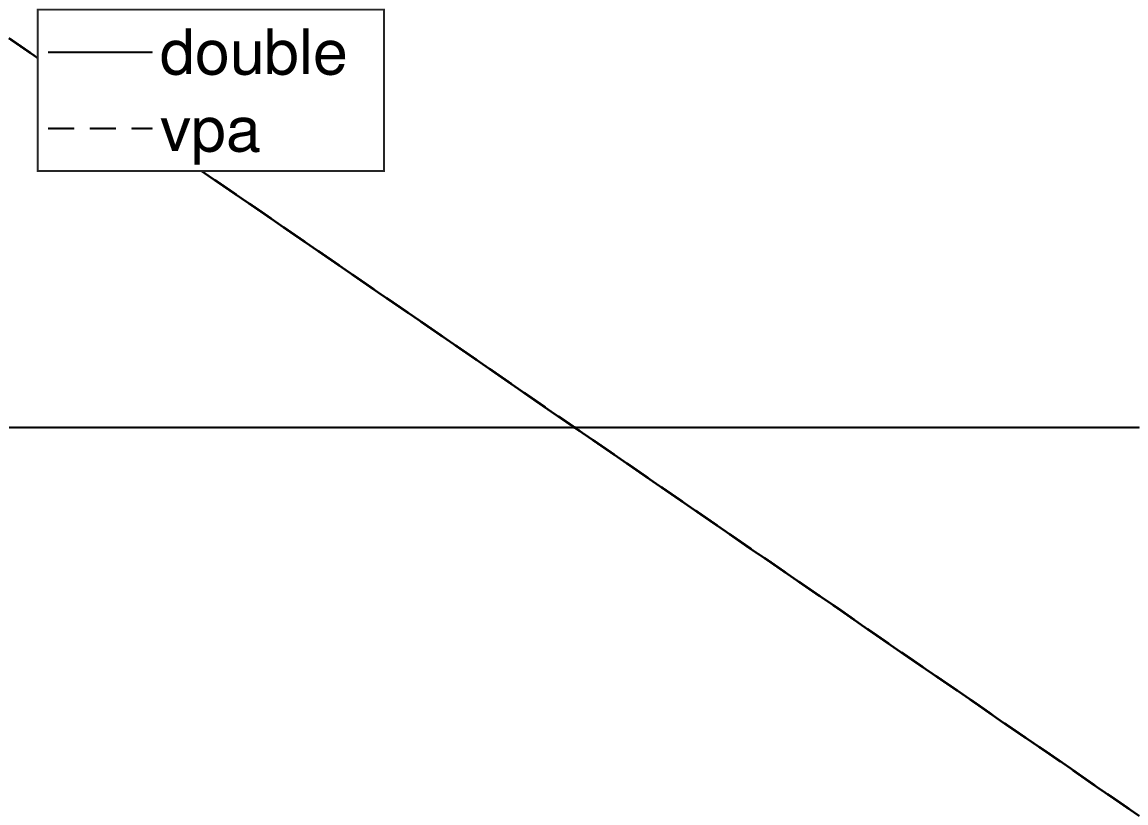}
			\caption{Eigenfrequency $\omega_1$}
			\label{fig:frequ1}
		\end{minipage}
		\hfil
		\begin{minipage}{.4\textwidth}
			\centering
			\includegraphics[scale=0.3]{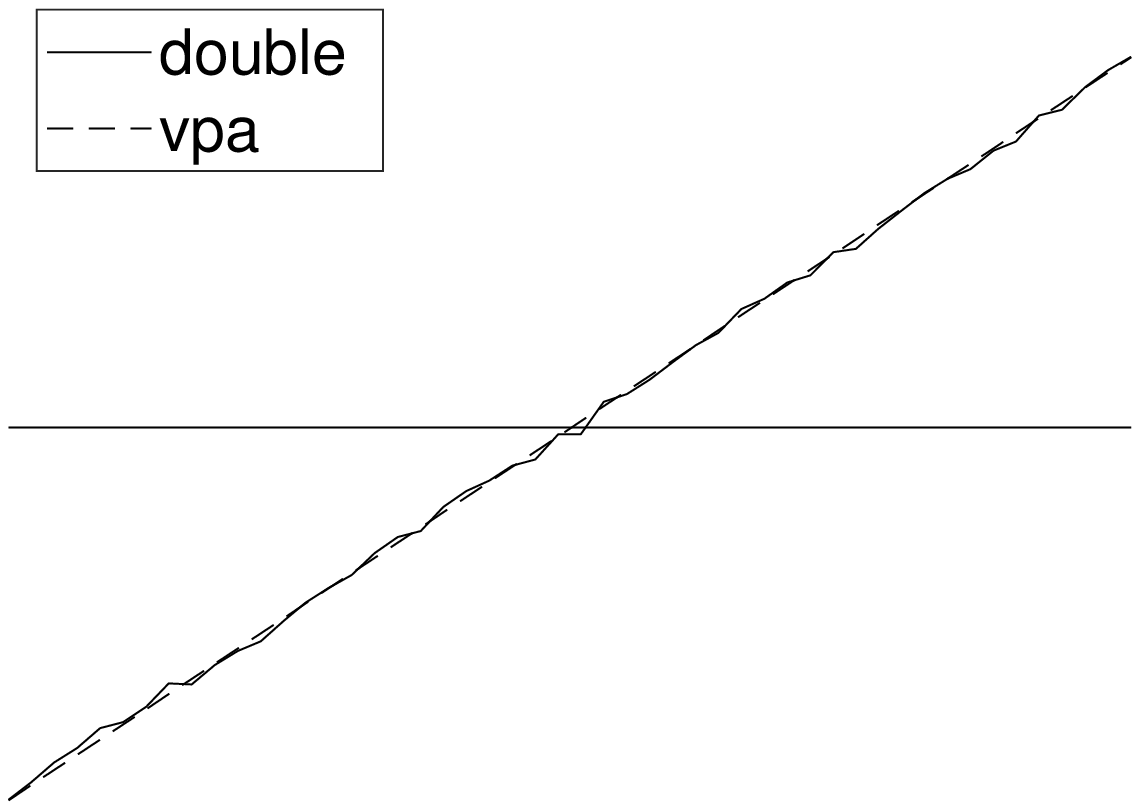}
			\caption{Eigenfrequency $\omega_2$}
			\label{fig:frequ2}
		\end{minipage}
		
		\begin{minipage}{.4\textwidth}
			\centering
			\includegraphics[scale=0.3]{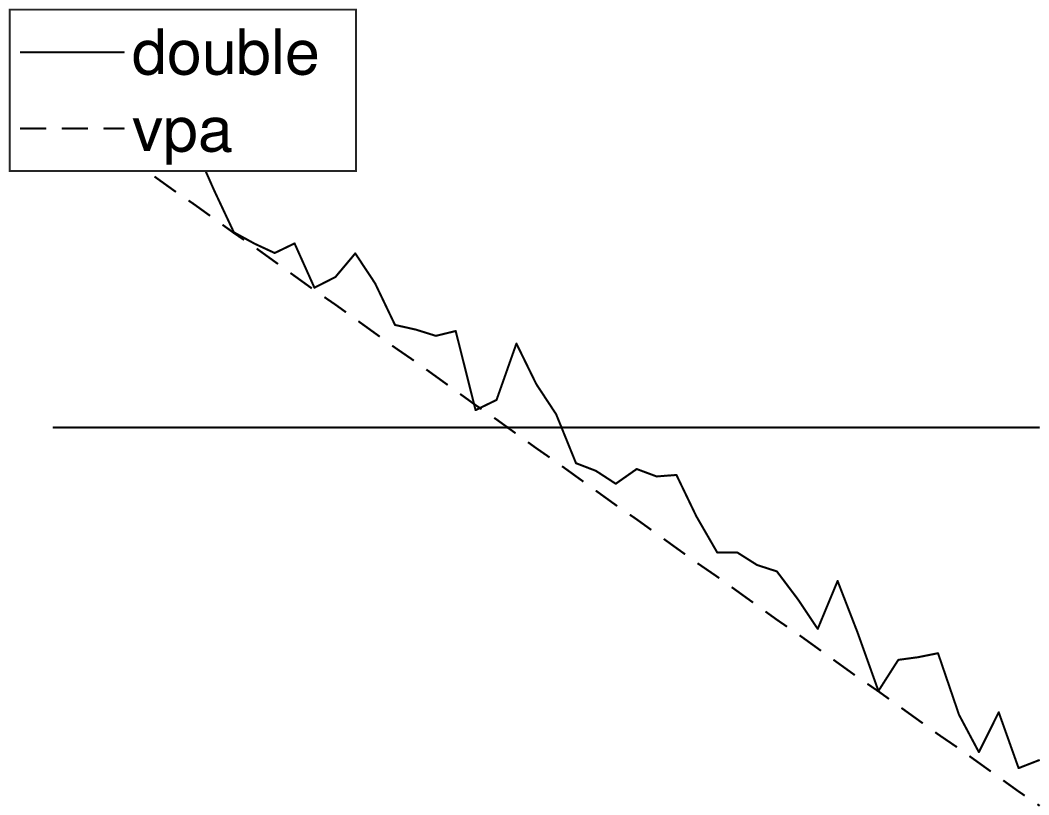}
			\caption{Eigenfrequency $\omega_3$}
			\label{fig:frequ3}
		\end{minipage}
		\hfil
		\begin{minipage}{.4\textwidth}
			\centering
			\includegraphics[scale=0.3]{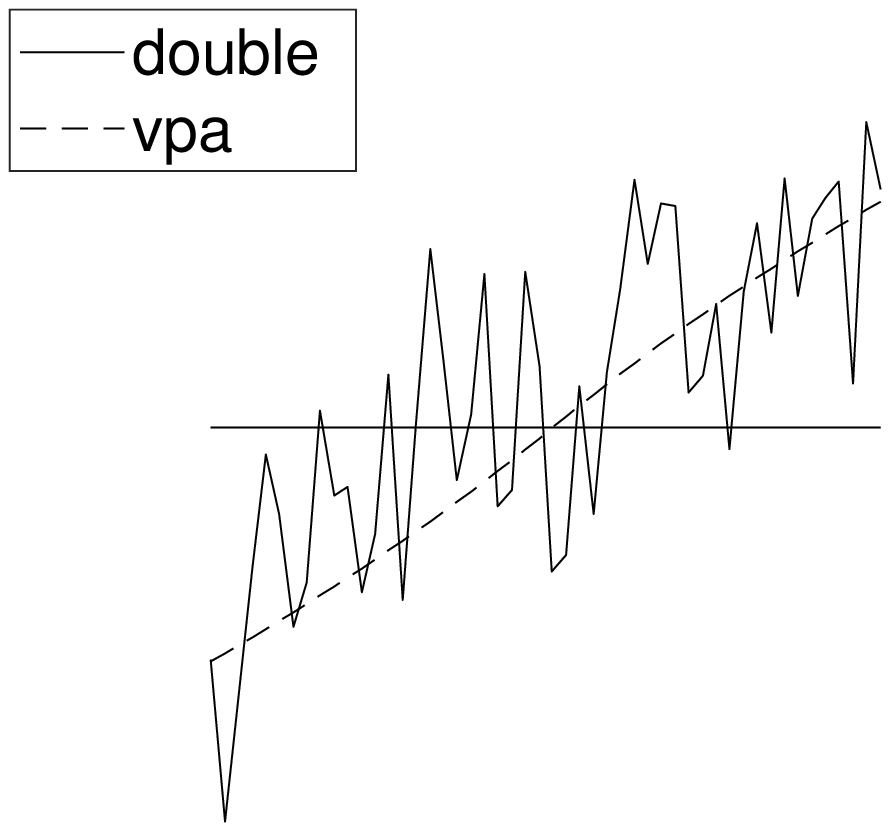}
			\caption{Eigenfrequency $\omega_4$ }
			\label{fig:frequ4}
		\end{minipage}
	\end{subfigure}
	\caption{Zoomed-in pictures of the frequency function $f(\omega)$ at the eigenfrequencies $\omega_i$. In all graphs $f(x_i)$ is ploted with $x_i \in [\omega_i-\epsilon; \ \omega_i+\epsilon]$, $\epsilon =10^{-12}$.
	}
	\label{fig:f2}
\end{figure}

When using the common trigonometric-hyperbolic ansatz \cref{for:eigenfunction} for arbitrary frame structures, one will sooner obtain inadequate results due to numerical inaccuracy when evaluating the frequency function $f(\omega)$. Especially, when the wanted natural frequency lies in a rather high-frequency range and/or the frequency function becomes lengthy e.g. for more complicated structures. Thus, numerical precision plays an important role. However, the computational performance can be improved by referring to the \emph{alternative ansatz} \cref{for:eigenfunction2}. With this ansatz, the numerical results are more stable for higher frequencies and rather complicated frequency functions and do not rely on a very high numerical precision. To illustrate this, in \Cref{fig:f2}, zoomed-in pictures show the frequency function $f(\omega)$ based on ansatz \cref{for:eigenfunction} of a frame structure in the area of the first four natural frequencies. Here, the frequency function is once evaluated with double precision and once using the variable precision arithmetic $vpa$ in Matlab. We observe that for higher frequencies the frequency function becomes rather jagged for double precision. As a result, one can obtain multiple natural frequencies where there is only one expected.
\subsection{Frame Structure}
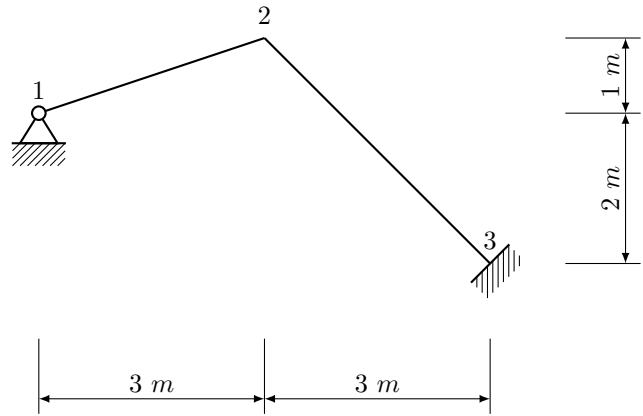
\begin{figure}[H]	
		\centering
		\begin{tikzpicture}[scale=1]
\draw[thick] (0.00000,0.00000) -- (3.00000,1.00000);
\draw[thick] (3.00000,1.00000) -- (6.00000,-2.00000);
\node at (0,0.3){1};
\node at (3,1.3){2};
\node at (6,-1.7){3};
\Festlager {(0.00000,0.00000)}[0.000000][1];
\Einspannung {(6.00000,-2.00000)}[45][1];
%

\draw(0,-3)--(0,-4);
\draw(3,-3)--(3,-4);
\draw(6,-3)--(6,-4);
\draw[latex-latex](0,-3.8)--(3,-3.8);
\draw[latex-latex](3,-3.8)--(6,-3.8);

\draw(7,-2)--(8,-2);
\draw(7,0)--(8,0);
\draw(7,1)--(8,1);
\draw[latex-latex](7.8,-2)--(7.8,0);
\draw[latex-latex](7.8,1)--(7.8,0);

\node[rotate=90] at (7.6,0.5){ $1 \ m$};
\node[rotate=90] at (7.6,-1){ $2 \ m$};

\node[rotate=0] at (1.5,-3.6){ $3 \ m$};
\node[rotate=0] at (4.5,-3.6){ $3 \ m$};

\end{tikzpicture}
		\caption{The modeled frame structure.}
		\label{fig:twoBeam}		
\end{figure}

We present a test example of a structure consisting of two rigidly connected beams forming an arbitrary frame structure. The frame is pinned on one side and clamped on the other side, see \Cref{fig:twoBeam}. 
The test example was modeled based on a steel profile, with the parameters $ A_m= 7.56 \cdot 10^{-4} \ [m^2]$, $ I_m = 3.5 \cdot 10^{-10} \ [m^4] $, $ \rho_m =  7.85 \ [t/m^3]$ and $ E_m = 2 \cdot 10^{8} \ [kN/m^2]$ for each beam.
The solution procedure was applied to this structure by using the three ansatz functions (a)$\rightarrow$\cref{for:eigenfunction}, (b)$\to$\cref{for:eigenfunction2.1}, (c)$\to$\cref{for:eigenfunction2} and the results were compared to a finite element computation with an alternating size of elements per beam, see \Cref{tab:naturalFrequTwoBeamFrame}. Also, the times to compute natural frequencies were compared between the NAT and the FEM approach, see \Cref{tab:timesTwoBeamFrame}. There, the computational time was taken for the last natural frequency computed using the FE Method that still lies within a range of a 10\% deviation from the NAT result with ansatz (c) and compared to the computation time using the NAT with ansatz (c). Further, the condition number was plotted for this example and can be studied in \Cref{fig:condNumTBF2}.
\begin{figure}[H]
	\centering
	\scalebox{0.6}{\hspace*{-8em}
	\begin{tikzpicture}[scale=1.7,shift={(-1em,0)}]
		\draw[dashed,draw=gray] (3.00000,1.00000) -- (0.00000,0.00000);
		\draw plot [smooth,thick] coordinates{ (-0.00000,0.00000)  (0.03551,-0.00551)  (0.07101,-0.01100)  (0.10650,-0.01646)  (0.14197,-0.02188)  (0.17743,-0.02724)  (0.21286,-0.03252)  (0.24826,-0.03772)  (0.28363,-0.04281)  (0.31896,-0.04779)  (0.35424,-0.05263)  (0.38948,-0.05733)  (0.42466,-0.06187)  (0.45979,-0.06623)  (0.49485,-0.07041)  (0.52985,-0.07439)  (0.56477,-0.07815)  (0.59962,-0.08168)  (0.63439,-0.08498)  (0.66907,-0.08802)  (0.70366,-0.09079)  (0.73817,-0.09328)  (0.77257,-0.09549)  (0.80688,-0.09739)  (0.84107,-0.09898)  (0.87517,-0.10024)  (0.90914,-0.10117)  (0.94301,-0.10174)  (0.97675,-0.10197)  (1.01037,-0.10182)  (1.04387,-0.10130)  (1.07723,-0.10039)  (1.11047,-0.09908)  (1.14357,-0.09737)  (1.17653,-0.09525)  (1.20935,-0.09270)  (1.24203,-0.08973)  (1.27457,-0.08633)  (1.30695,-0.08248)  (1.33919,-0.07818)  (1.37128,-0.07343)  (1.40321,-0.06821)  (1.43499,-0.06253)  (1.46661,-0.05639)  (1.49807,-0.04976)  (1.52937,-0.04266)  (1.56051,-0.03507)  (1.59149,-0.02700)  (1.62231,-0.01844)  (1.65297,-0.00940)  (1.68345,0.00015)  (1.71378,0.01018)  (1.74394,0.02071)  (1.77393,0.03174)  (1.80376,0.04326)  (1.83343,0.05527)  (1.86293,0.06778)  (1.89227,0.08079)  (1.92144,0.09428)  (1.95045,0.10827)  (1.97929,0.12274)  (2.00798,0.13769)  (2.03650,0.15313)  (2.06487,0.16905)  (2.09307,0.18543)  (2.12112,0.20229)  (2.14902,0.21961)  (2.17676,0.23739)  (2.20436,0.25563)  (2.23180,0.27431)  (2.25909,0.29343)  (2.28625,0.31299)  (2.31326,0.33297)  (2.34012,0.35337)  (2.36686,0.37418)  (2.39345,0.39540)  (2.41992,0.41702)  (2.44626,0.43902)  (2.47247,0.46139)  (2.49856,0.48414)  (2.52452,0.50724)  (2.55038,0.53070)  (2.57612,0.55449)  (2.60175,0.57861)  (2.62727,0.60305)  (2.65269,0.62779)  (2.67801,0.65284)  (2.70324,0.67816)  (2.72838,0.70376)  (2.75343,0.72962)  (2.77839,0.75574)  (2.80328,0.78209)  (2.82809,0.80866)  (2.85283,0.83545)  (2.87750,0.86245)  (2.90211,0.88963)  (2.92666,0.91699)  (2.95116,0.94451)  (2.97560,0.97219)  (3.00000,1.00000) };
		\Festlager {(0.00000,-0.00000)}[0.000000][1.000000];
		\draw[dashed,draw=gray] (3.00000,1.00000) -- (6.00000,-2.00000);
		\draw plot [smooth] coordinates{ (6.00000,-2.00000)  (5.97025,-1.96914)  (5.94159,-1.93720)  (5.91396,-1.90422)  (5.88732,-1.87026)  (5.86162,-1.83535)  (5.83682,-1.79954)  (5.81288,-1.76288)  (5.78974,-1.72541)  (5.76736,-1.68718)  (5.74570,-1.64824)  (5.72472,-1.60862)  (5.70435,-1.56838)  (5.68457,-1.52755)  (5.66532,-1.48619)  (5.64657,-1.44434)  (5.62826,-1.40205)  (5.61035,-1.35935)  (5.59280,-1.31629)  (5.57557,-1.27291)  (5.55861,-1.22927)  (5.54188,-1.18539)  (5.52534,-1.14133)  (5.50894,-1.09712)  (5.49265,-1.05280)  (5.47642,-1.00843)  (5.46022,-0.96402)  (5.44400,-0.91964)  (5.42773,-0.87530)  (5.41136,-0.83106)  (5.39487,-0.78695)  (5.37821,-0.74300)  (5.36135,-0.69926)  (5.34425,-0.65575)  (5.32689,-0.61251)  (5.30922,-0.56957)  (5.29122,-0.52696)  (5.27285,-0.48472)  (5.25409,-0.44288)  (5.23491,-0.40146)  (5.21527,-0.36048)  (5.19516,-0.31999)  (5.17455,-0.28000)  (5.15341,-0.24053)  (5.13172,-0.20161)  (5.10947,-0.16326)  (5.08662,-0.12550)  (5.06316,-0.08835)  (5.03908,-0.05183)  (5.01435,-0.01595)  (4.98897,0.01927)  (4.96291,0.05382)  (4.93617,0.08768)  (4.90874,0.12086)  (4.88060,0.15333)  (4.85175,0.18508)  (4.82219,0.21613)  (4.79190,0.24644)  (4.76089,0.27604)  (4.72914,0.30490)  (4.69667,0.33304)  (4.66348,0.36045)  (4.62955,0.38713)  (4.59491,0.41309)  (4.55955,0.43833)  (4.52347,0.46287)  (4.48670,0.48670)  (4.44923,0.50983)  (4.41107,0.53229)  (4.37225,0.55407)  (4.33277,0.57519)  (4.29264,0.59567)  (4.25189,0.61552)  (4.21052,0.63476)  (4.16857,0.65341)  (4.12603,0.67149)  (4.08295,0.68901)  (4.03933,0.70600)  (3.99520,0.72247)  (3.95059,0.73847)  (3.90552,0.75400)  (3.86001,0.76910)  (3.81409,0.78378)  (3.76779,0.79809)  (3.72113,0.81204)  (3.67415,0.82567)  (3.62688,0.83900)  (3.57934,0.85207)  (3.53157,0.86490)  (3.48359,0.87753)  (3.43545,0.88999)  (3.38717,0.90232)  (3.33878,0.91454)  (3.29032,0.92668)  (3.24183,0.93880)  (3.19333,0.95090)  (3.14486,0.96304)  (3.09646,0.97525)  (3.04817,0.98756)  (3.00000,1.00000) };
		\Einspannung {(6.00000,-2.00000)}[45.000000][1.000000];
	\end{tikzpicture}}
	\caption{The first mode shape of the frame structure.}
	\label{fig:modeShapeTBF}
\end{figure}
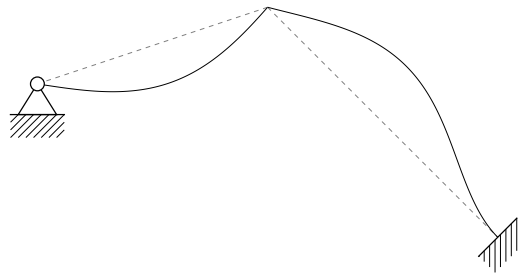
\begin{figure}[H]
	\centering
	\scalebox{0.6}{\input{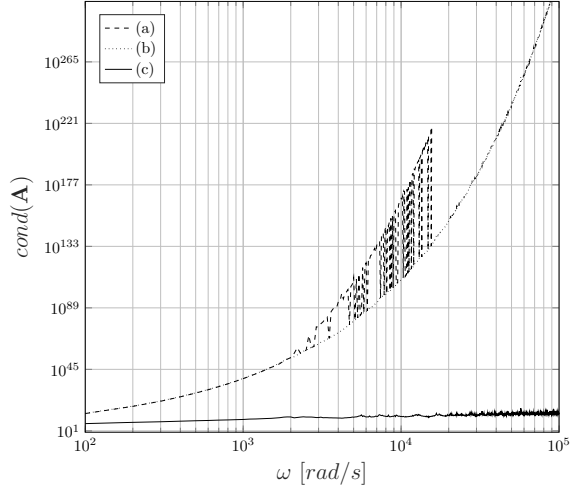}}
	\caption{The condition number of the system matrix of the frame structure.}
	\label{fig:condNumTBF2}
\end{figure}

\begin{table}[H]
\centering
\caption{The computation times of the frame structure.}
\label{tab:timesTwoBeamFrame}
\begin{tabular}{|c|c|c|c|c|}\hline
 \multirow{2}{*}{\shortstack{Number of\\ Frequency}} & \multicolumn{3}{c|}{FEM} &NAT\\\cline{2-5}
& max. $\omega [rad/s]$ \tablefootnote{last frequency computed with FEM that has less than 10\% deviation from NAT computation with Ansatz (c)}& Elements per beam & t[s] & t[s] \\\hline
2& 4.8627& 2& $<0.1$& $<0.1$ \\6& 33.0061& 4& $<0.1$& $<0.1$ \\13& 128.3669& 8& $<0.1$& $<0.1$ \\32& 716.3521& 16& $<0.1$& $<0.1$ \\73& 3510.2522& 32& $<0.1$& 0.12569 \\151& 14234.4738& 64& 0.1068& 0.43859 \\320& 59048.7810& 128& 0.9619& 0.51866 \\702& 243982.2853& 256& 15.7096& 2.2451 \\1594& 957605.9388& 512& 192.2587& 42.2335 \\3195& 2771277.8678& 1024& 1595.0565& 42.2877 \\\hline
\end{tabular}
\end{table}

\begin{table}[H]
\centering
\caption{The natural frequencies of the frame structure in [rad/s].}
\label{tab:naturalFrequTwoBeamFrame}
\setlength\extrarowheight{5pt}
\resizebox{\columnwidth}{!}{%
\begin{tabular}{|c|c|c|c|c|c|c|c|c|c|}\hline
\multirow{3}{*}{\shortstack{Number of\\ Natural\\ Frequency}}& \multicolumn{6}{c|}{NAT} & \multicolumn{3}{c|}{FEM}\\\cline{2-10}
& \multicolumn{3}{c|}{double} & \multicolumn{3}{c|}{vpa} & \multicolumn{3}{c|}{Elements per Beam} \\\cline{2-10}
& (a)  & (b) & (c) & (a)  & (b) & (c) & 16&32&64 \\ \hline

1& 3.1094& 3.1094& 3.1094& 3.1094& 3.1094& 3.1094& 3.1094& 3.1094& 3.1094\\
2& 4.8078& 4.8078& 4.8078& 4.8078& 4.8078& 4.8078& 4.8078& 4.8078& 4.8078\\
3& 10.4144& 10.4144& 10.4144& 10.4144& 10.4144& 10.4144& 10.414\underline{7}& 10.4144& 10.4144\\
&..&..&..&..&..&..&..&..&..\\
15& 162.81\underline{58}& 162.8160& 162.8160& 162.8160& 162.8160& 162.8160& 16\underline{3.7454}& 162.8\underline{787}& 162.8\underline{200}\\
16& 174.603\underline{0}& 174.6036& 174.6036& 174.6036& 174.6036& 174.6036& 17\underline{5.3495}& 174.6\underline{525}& 174.60\underline{67}\\
17& 20\underline{1.9954}& 202.0324& 202.0324& 202.0324& 202.0324& 202.0324& 20\underline{3.9897}& 202.\underline{1665}& 202.0\underline{409}\\\cline{2-2}
&&..&..&..&..&..&..&..&..\\
31 && 634.9490& 634.9490& 634.9490& 634.9490& 634.9490& 6\underline{88.0855}& 63\underline{8.8745}& 63\underline{5.2119}\\
32 && 675.0622& 675.0622& 675.0622& 675.0622& 675.0622& \underline{716.3521}& 67\underline{7.3190}& 675.\underline{2117}\\
33 && 708.8620& 708.8620& 708.8620& 708.8620& 708.8620& \pur{787.7197}& 7\underline{14.2016}& 70\underline{9.2205}\\
&&..&..&..&..&..&..&..&..\\
48 && 1498.5735& 1498.5735& 1498.\underline{3202}& 1498.5735& 1498.5735& \pur{1848.5917}& 1\underline{523.7303}& 1\underline{500.7542}\\
49 && 1526.1712& 1526.1712& 152\underline{7.1223}& 1526.1712& 1526.1712& \pur{1923.3839}& 15\underline{65.8748}& 152\underline{8.9541}\\
50 && 1618.8520& 1618.8520& 1\underline{581.1209}& 1618.8520& 1618.8520& \pur{2003.3543}& 16\underline{63.6016}& 16\underline{22.9055}\\\cline{5-5}
&&&..&&..&&..&..&..\\
72 && 3099.1823& 3099.1823 && 3099.1823& 3099.1823& \pur{21006.9000}& 3\underline{386.8400}& 3\underline{120.1627}\\
73 && 3219.5726& 3219.5726 && 3219.5726& 3219.5726& \pur{22796.4503}& 3\underline{510.2522}& 32\underline{48.4805}\\
74 && 3289.0796& 3289.0796 && 3289.0796& 3289.0796& \pur{25373.8867}& \pur{3711.2534}& 3\underline{306.7830}\\\cline{8-8}
&&..&..&&..&..&&..&..\\
150 && 5146.4523& 5146.4523 && 5146.4523& 5146.4523& & \pur{53266.5249}& \underline{13871.6132}\\
151 && 5198.0249& 5198.0249 && 5198.0249& 5198.0249& & \pur{54623.1788}& \underline{14234.4738}\\
152 && 5357.7191& 5357.7191 && 5357.7191& 5357.7191& & \pur{57991.3441}& \pur{14646.4119}\\
&&..&..&&..&..&&..&..\\
188 && 20082.8678& 20082.8678 && 20082.8678& 20082.8678 && \pur{171783.6588}& \pur{23671.7495}\\
189 && 20161.3979& 20161.3979 && 20161.3979& 20161.3979 && \pur{174609.8040}& \pur{23903.3555}\\
190 && 20449.5667& 20449.5667 && 20449.5667& 20449.5667 && \pur{176350.0940}& \pur{24275.7074}\\\cline{9-9}
&&..&..&&..&..&\multicolumn{2}{c|}{}&..\\
191 && 20516.9051& 20516.9051 && 20516.9051& 20516.9051 & \multicolumn{2}{c|}{}& \pur{24345.0237}\\
192 && 20733.0092& 20733.0092 && 20733.0092& 20733.0092 & \multicolumn{2}{c|}{}& \pur{24662.7192}\\
193 && 20919.7970& 20919.7970 && 20919.7970& 20919.7970 & \multicolumn{2}{c|}{}& \pur{25015.2626}\\\cline{3-3}
&\multicolumn{2}{c|}{}&..&&..&..&\multicolumn{2}{c|}{}&..\\
219 & \multicolumn{2}{c|}{}& 26794.1655 && 26794.1655& 26794.1655 & \multicolumn{2}{c|}{}& \pur{33423.5651}\\
220 & \multicolumn{2}{c|}{}& 26915.9936 && 26915.9936& 26915.9936 & \multicolumn{2}{c|}{}& \pur{33780.5983}\\
221 & \multicolumn{2}{c|}{}& 27220.7782 & &27220.7782 &27220.7782 & \multicolumn{2}{c|}{}& \pur{34152.4567}\\\cline{6-6}
&\multicolumn{2}{c|}{}&..&\multicolumn{2}{c|}{}&..&\multicolumn{2}{c|}{}&..\\
380 & \multicolumn{2}{c|}{}& 74088.3050 & \multicolumn{2}{c|}{} & 74088.3050& \multicolumn{2}{c|}{}& \pur{351135.7772}\\
381 & \multicolumn{2}{c|}{}& 74411.2521 & \multicolumn{2}{c|}{} & 74411.2521& \multicolumn{2}{c|}{}& \pur{352652.5927}\\
382 & \multicolumn{2}{c|}{}& 74862.5654 & \multicolumn{2}{c|}{} & 74862.5654& \multicolumn{2}{c|}{}& \pur{353569.1257}\\\cline{10-10}
&\multicolumn{2}{c|}{}&..&\multicolumn{2}{c|}{}&..&\multicolumn{3}{c|}{}\\
1735 & \multicolumn{2}{c|}{}& 997062.8183 & \multicolumn{2}{c|}{} & 997062.8183& \multicolumn{3}{c|}{}\\
1736 & \multicolumn{2}{c|}{}& 998652.9269 & \multicolumn{2}{c|}{} & 998652.9269& \multicolumn{3}{c|}{}\\
1737 & \multicolumn{2}{c|}{}& 999016.8478 & \multicolumn{2}{c|}{} & 999016.8478& \multicolumn{3}{c|}{}\\
&\multicolumn{2}{c|}{}&..&\multicolumn{2}{c|}{}&..&\multicolumn{3}{c|}{}\\
\hline
\end{tabular}}
\end{table}

In \Cref{fig:condNumTBF2} it can be observed that the condition number for the system matrix of the frame structure using ansatz (c) stays low even for high frequencies which is an indicator for a stable system of equations. Whereas the condition number for the system matrices using ansatz (a) and (b) increases drastically with rising frequency $\omega$. The commonly used ansatz (a) even forms additional peaks before reaching complete failure. 
\Cref{tab:naturalFrequTwoBeamFrame} shows that the natural frequencies computed with the NAT are more accurate than the FEM results. Numbers are underlined if they differ from the NAT solution using ansatz (c). A result is colored gray when it passed the 10\% deviation from the NAT solution using ansatz (c). When there is a blank cell in the table, no result was able to be computed, i.e. the computation failed. It is also noticeable that regarding the NAT, the ansatz (c), the presented \textit{alternative ansatz}, performs better than the other two general solutions. While even using a 'lower' double-precision, the NAT with ansatz (c) outperforms the FEM as well as the other two ansatz functions. As noticed, the two dots in the last row of the columns for ansatz (c) indicate that the limit has not yet been reached there, so even higher natural frequencies can be computed using this ansatz. 
In \Cref{tab:timesTwoBeamFrame} the computational times for the FEM look promising for lower natural frequencies but rise tremendously as the frequency range increases. Contrarily, the times to compute the same natural frequencies using the NAT with ansatz (c) increase rather moderately and stay low for immensely high frequencies.
\subsection{Bridge}
A reinforced concrete bridge was modeled regarding the parameters $L=20 [m]$, $E=31 \cdot 10^6 [kN/m^2]$, $\rho=2.3 [t/m^3]$ for all beams. The beams \fbox{\begin{minipage}{.4em}1\end{minipage}} - \fbox{\begin{minipage}{.4em}3\end{minipage}} share the parameters $A=12 [m^2]$, $I=40 [m^4]$ and for the remaining beams \fbox{\begin{minipage}{.4em}4\end{minipage}} - \fbox{\begin{minipage}{.4em}5\end{minipage}} the parameters are $A=9 [m^2]$, $I=30 [m^4]$.
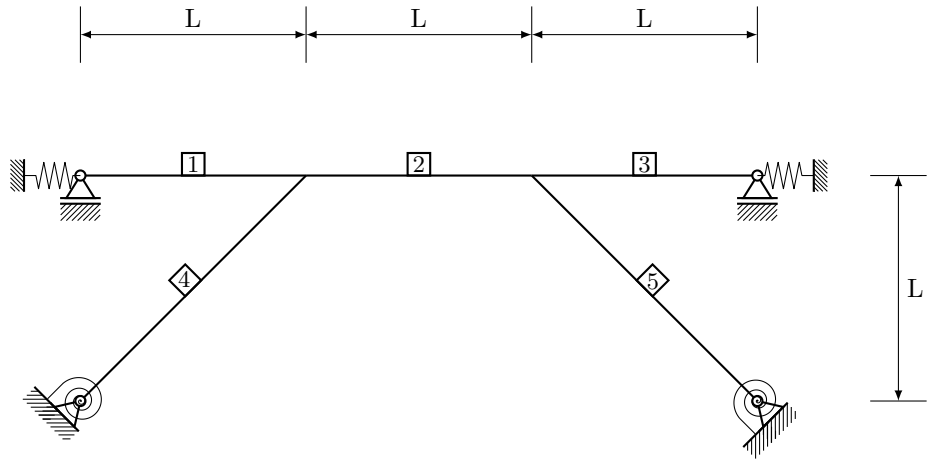
\begin{figure}[H]
	\centering
	\begin{tikzpicture}[scale=0.15]
	\draw[thick] (9,0)rectangle(11,2);
	\node at (10,1){\small1};
	\draw[thick] (29,0)rectangle(31,2);
	\node at (30,1){\small2};
	\draw[thick] (49,0)rectangle(51,2);
	\node at (50,1){\small3};
	\draw[thick,rotate around={45:(0,-20)}] (28.28/2-1,-20)rectangle(28.28/2+1,-18);
	\node at (9.2,-9.2){\small4};
	\draw[thick,rotate around={-45:(60,-20)}] (60-28.28/2-1,-20)rectangle(60-28.28/2+1,-18);
	\node at (40+10.75,-9.4){\small5};
	\draw[thick] (0.00000,0.00000) -- (20.00000,0.00000);
	\draw[thick] (20.00000,0.00000) -- (40.00000,0.00000);
	\draw[thick] (40.00000,0.00000) -- (60.00000,0.00000);
	\draw[thick] (0.00000,-20.00000) -- (20.00000,0.00000);
	\draw[thick] (60.00000,-20.00000) -- (40.00000,0.00000);
	\Loslager{(0.00000,0.00000)}[0.000000][5.000000];
	\Loslager {(60.00000,0.00000)}[0.000000][5.000000];
	\Festlager {(0.00000,-20.00000)}[-45.000000][5.000000];
	\Festlager {(60.00000,-20.00000)}[45.000000][5.000000];
	\begin{scope}[shift={(0,0)},rotate=180,scale=2]
	\draw[snake=zigzag,line before snake=1mm, line after snake=1mm,
	segment length=4pt,
	segment amplitude=5pt,join=round, thin] 
	(0,0.) -- (2.5,0.); 
	\Einspannung{(2.5,0)}[90][2][180]
	{}
	\end{scope}
	\begin{scope}[shift={(60,0)},rotate=0,scale=2]
	\draw[snake=zigzag,line before snake=1mm, line after snake=1mm,
	segment length=4pt,
	segment amplitude=5pt,join=round, thin] 
	(0,0.) -- (2.5,0.); 
	\Einspannung{(2.5,0)}[90][2][0]%
	{}
	\end{scope}
	\Drehfeder {(0.00000,-20.00000)}[-45.000000][5.000000];
	\Drehfeder {(60.00000,-20.00000)}[45.000000][5.000000];
	\draw (0,10)--(0,15);
	\draw (20,10)--(20,15);
	\draw (40,10)--(40,15);
	\draw (60,10)--(60,15);
	\draw[latex-latex](0,12.5)--(20,12.5);
	\draw[latex-latex](20,12.5)--(40,12.5);
	\draw[latex-latex](40,12.5)--(60,12.5);
	\node at (10,14){L};
	\node at (30,14){L};
	\node at (50,14){L};
	\draw(70,0)--(75,0);
	\draw(70,-20)--(75,-20);
	\draw[latex-latex](72.5,0)--(72.5,-20);
	\node at (74,-10){L};
\end{tikzpicture}
	\caption{The modeled bridge.}
	\label{fig:bridgeSketch}
\end{figure}
\begin{figure}[H]
	\centering
	\begin{tikzpicture}[scale=0.15]
\draw[dotted,draw=gray] (0.00000,0.00000) -- (20.00000,0.00000);
\draw plot [smooth] coordinates{ (-3.17414,-0.00000)  (-2.12133,0.22776)  (-1.06813,0.45349)  (-0.01456,0.67515)  (1.03938,0.89074)  (2.09370,1.09828)  (3.14839,1.29581)  (4.20346,1.48144)  (5.25890,1.65331)  (6.31472,1.80962)  (7.37091,1.94864)  (8.42747,2.06872)  (9.48440,2.16826)  (10.54171,2.24578)  (11.59938,2.29984)  (12.65743,2.32913)  (13.71585,2.33239)  (14.77463,2.30849)  (15.83379,2.25637)  (16.89331,2.17506) };
\Loslager {(-3.17414,-0.00000)}[0.000000][5.000000];
\begin{scope}[shift={(-3.17414,0)},rotate=180,scale=2]
	\draw[snake=zigzag,line before snake=1mm, line after snake=1mm,
	segment length=2pt,
	segment amplitude=6pt,join=round, thin] 
	(0,0.) -- (2,0.); 
	\Einspannung{(2,0)}[90][2][180]%
	{}
\end{scope}
\draw[dotted,draw=gray] (20.00000,0.00000) -- (40.00000,0.00000);
\draw plot [smooth] coordinates{ (16.89331,2.17506)  (17.94263,2.05887)  (18.99231,1.90546)  (20.04236,1.71899)  (21.09278,1.50372)  (22.14357,1.26394)  (23.19472,1.00403)  (24.24625,0.72842)  (25.29814,0.44156)  (26.35040,0.14794)  (27.40304,-0.14794)  (28.45604,-0.44156)  (29.50941,-0.72842)  (30.56314,-1.00403)  (31.61725,-1.26394)  (32.67173,-1.50372)  (33.72657,-1.71899)  (34.78178,-1.90546)  (35.83737,-2.05887)  (36.89331,-2.17506) };
\draw[dotted,draw=gray] (40.00000,0.00000) -- (60.00000,0.00000);
\draw plot [smooth] coordinates{ (36.89331,-2.17506)  (37.93905,-2.25637)  (38.98516,-2.30849)  (40.03164,-2.33239)  (41.07848,-2.32913)  (42.12570,-2.29984)  (43.17328,-2.24578)  (44.22124,-2.16826)  (45.26957,-2.06872)  (46.31827,-1.94864)  (47.36735,-1.80962)  (48.41680,-1.65331)  (49.46662,-1.48144)  (50.51681,-1.29581)  (51.56738,-1.09828)  (52.61833,-0.89074)  (53.66965,-0.67515)  (54.72134,-0.45349)  (55.77341,-0.22776)  (56.82586,-0.00000) };
\Loslager {(56.82586,-0.00000)}[0.000000][5.000000];
\begin{scope}[shift={(56.82586,0)},rotate=0,scale=2]
	\draw[snake=zigzag,line before snake=1mm, line after snake=1mm,
	segment length=6pt,
	segment amplitude=3pt,join=round, thin] 
	(0,0.) -- (3.5,0.); 
	\Einspannung{(3.5,0)}[90][2][0]%
	{}
\end{scope}
\draw[dotted,draw=gray] (0.00000,-20.00000) -- (20.00000,0.00000);
\draw plot [smooth] coordinates{ (-0.00000,-20.00000)  (0.72986,-18.67433)  (1.46313,-17.35206)  (2.20318,-16.03655)  (2.95329,-14.73106)  (3.71658,-13.43871)  (4.49602,-12.16245)  (5.29434,-10.90499)  (6.11400,-9.66879)  (6.95717,-8.45601)  (7.82570,-7.26849)  (8.72107,-6.10768)  (9.64437,-4.97468)  (10.59630,-3.87017)  (11.57711,-2.79439)  (12.58666,-1.74718)  (13.62432,-0.72791)  (14.68905,0.26448)  (15.77936,1.23149)  (16.89331,2.17506) };
\Festlager {(-0.00000,-20.00000)}[-45.000000][5.000000];
\Drehfeder {(-0.00000,-20.00000)}[-45.000000][5.000000];
\draw[dotted,draw=gray] (60.00000,-20.00000) -- (40.00000,0.00000);
\draw plot [smooth] coordinates{ (60.00000,-20.00000)  (58.62460,-19.22041)  (57.25260,-18.43741)  (55.88739,-17.64766)  (54.53224,-16.84789)  (53.19027,-16.03497)  (51.86444,-15.20598)  (50.55750,-14.35817)  (49.27189,-13.48911)  (48.00981,-12.59662)  (46.77307,-11.67888)  (45.56318,-10.73442)  (44.38122,-9.76216)  (43.22788,-8.76141)  (42.10343,-7.73193)  (41.00771,-6.67388)  (39.94011,-5.58788)  (38.89958,-4.47500)  (37.88462,-3.33675)  (36.89331,-2.17506) };
\Festlager {(60.00000,-20.00000)}[45.000000][5.000000];
\Drehfeder {(60.00000,-20.00000)}[45.000000][5.000000];
\end{tikzpicture}
	\caption{The first mode shape of the bridge.}
	\label{fig:bridgeEigenmode}
\end{figure}
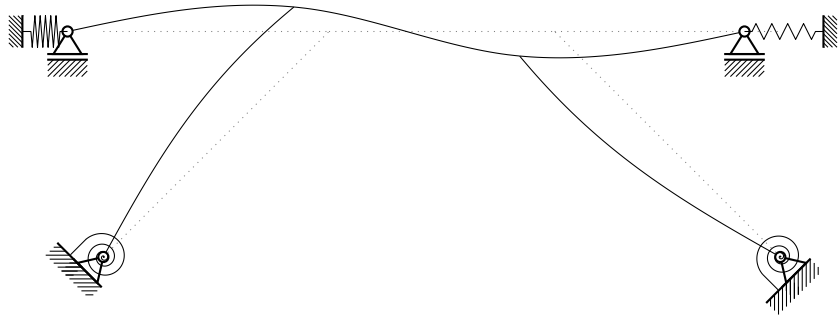
\begin{figure}[H]
	\centering
	\scalebox{0.6}{\input{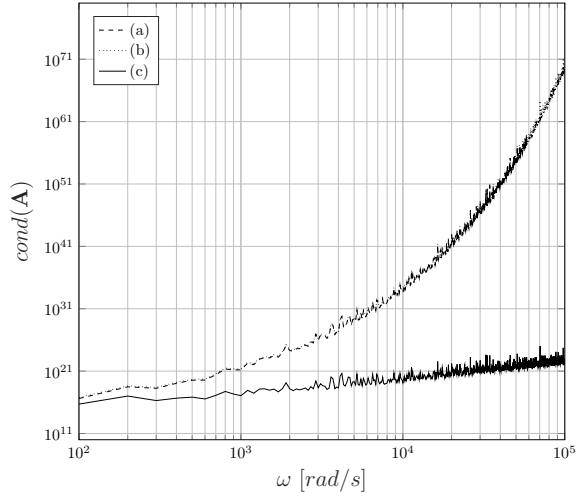}}
	\caption{The condition number of the system matrix of the bridge.}
	\label{fig:condNumBridge2}
\end{figure}

\begin{table}[H]
\centering
\caption{The computation times of the bridge.}
\label{tab:timesTwoBeamFrame}
\begin{tabular}{|c|c|c|c|c|}\hline
 \multirow{2}{*}{\shortstack{Number of\\ Frequency}} & \multicolumn{3}{c|}{FEM} &NAT\\\cline{2-5}
& max. $\omega [rad/s]$ \tablefootnote{last frequency computed with FEM that has less than 10\% deviation from NAT computation with Ansatz (c)}& Elements per beam & t[s] & t[s] \\\hline
11& 418.2554& 2& $<0.1$& $<0.1$ \\25& 1223.6150& 4& $<0.1$& $<0.1$ \\44& 2547.1030& 8& $<0.1$& 0.25898 \\74& 4837.9021& 16& $<0.1$& 1.0005 \\130& 9685.2435& 32& $<0.1$& 1.0574 \\223& 18003.4972& 64& $<0.1$& 4.5552 \\433& 38425.0172& 128& 1.0901& 4.7907 \\805& 75794.2360& 256& 12.9402& 22.3849 \\1586& 157012.4639& 512& 104.0063&    24.0491 \\3194& 330080.3245& 1024& 815.8332& 25.2141 \\\hline
\end{tabular}
\end{table}

\begin{table}[H]
\centering
\caption{The natural frequencies of the bridge in [rad/s].}
\label{tab:naturalFrequBridge2}
\setlength\extrarowheight{5pt}
\resizebox{\columnwidth}{!}{%
\begin{tabular}{|c|c|c|c|c|c|c|c|c|c|}\hline
\multirow{3}{*}{\shortstack{Number of\\Natural \\ Frequency}}& \multicolumn{6}{c|}{NAT} & \multicolumn{3}{c|}{FEM}\\\cline{2-10}
& \multicolumn{3}{c|}{double} & \multicolumn{3}{c|}{vpa} & \multicolumn{3}{c|}{Elements per Beam} \\\cline{2-10}
& (a)  & (b) & (c) & (a)  & (b) & (c) & 16&32&64 \\ \hline

1& 37.9854& 37.9854& 37.9854& 37.9854& 37.9854& 37.9854& 37.985\underline{6}& 37.9855& 37.9854\\
2& 74.2871& 74.2871& 74.2871& 74.2871& 74.2871& 74.2871& 74.287\underline{5}& 74.2872& 74.2871\\
3& 105.7450& 105.7450& 105.7450& 105.7450& 105.7450& 105.7450& 105.74\underline{87}& 105.745\underline{9}& 105.745\underline{3}\\
&..&..&..&..&..&..&..&..&..\\
73& 4461.7242& 4461.7242& 4461.7242& 4461.7242& 4461.7242& 4461.7242& 4\underline{753.2843}& 4\underline{539.9982}& 44\underline{80.5927}\\
74& 4473.8551& 4473.8551& 4473.8551& 4473.8551& 4473.8551& 4473.8551& 4\underline{837.9021}& 4\underline{634.3395}& 4\underline{515.8803}\\
75& 4515.0470& 4515.0470& 4515.0470& 4515.0470& 4515.0470& 4515.0470& \pur{5000.6678}& 4\underline{657.4950}& 45\underline{59.9003}\\
&..&..&..&..&..&..&..&..&..\\
129& 8812.4941& 8812.4941& 8812.4941& 8812.4941& 8812.4941& 8812.4941& \pur{12994.3482}& \underline{9646.0261}& 8\underline{965.7722}\\
130& 8853.3989& 8853.3989& 8853.3989& 8853.3989& 8853.3989& 8853.3989& \pur{13693.3704}& \underline{9685.2435}& \underline{9000.0630}\\
131& 8862.147\underline{3}& 8862.1475& 8862.1475& 8862.1475& 8862.1475& 8862.1475& \pur{14437.5183}& \pur{9878.7665}& \underline{9071.5678}\\\cline{2-2}
&&..&..&..&..&..&..&..&..\\
222 && 16416.1027& 16416.1027& 16416.1027& 16416.1027& 16416.1027& \pur{156994.8526}& \pur{20342.0473}& 1\underline{7874.0744}\\
223 && 16527.6906& 16527.6906& 16527.6906& 16527.6906& 16527.6906& \pur{165670.7114}& \pur{21219.2334}& 1\underline{8003.4972}\\
224 && 16545.6271& 16545.6271& 16545.6271& 16545.6271& 16545.6271& \pur{169638.2781}& \pur{21975.3448}& \pur{18225.9889}\\
&&..&..&..&..&..&..&..&..\\
235 && 17514.5070& 17514.5070& 17514.5070& 17514.5070& 17514.5070& \pur{212717.0285}& \pur{29380.2136}& \pur{19376.7279}\\
236 && 17545.6329& 17545.6329& 17545.6329& 17545.6329& 17545.6329& \pur{214599.6433}& \pur{30468.3455}& \pur{19426.6437}\\
237 && 17629.9634& 17629.9634& 17629.9634& 17629.9634& 17629.9634& \pur{214602.7984}& \pur{30985.1054}& \pur{19537.4222}\\\cline{8-8}
&&..&..&..&..&..&&..&..\\
475 && 38742.1566& 38742.1566& 38742.1566& 38742.1566& 38742.1566 && \pur{858548.1299}& \pur{122440.6046}\\
476 && 38762.7334& 38762.7334& 38762.7334& 38762.7334& 38762.7334 && \pur{860629.4273}& \pur{122704.7380}\\
477 && 38890.8164& 38890.8164& 38890.8164& 38890.8164& 38890.8164 && \pur{860630.2312}& \pur{123645.1451}\\\cline{9-9}
&&..&..&..&..&..&\multicolumn{2}{c|}{}&..\\
502 && 41198.3553& 41198.3553& 41198.3553& 41198.3553& 41198.3553 & \multicolumn{2}{c|}{}& \pur{166664.8265}\\
503 && 41202.1461& 41202.1461& 41202.1461& 41202.1461& 41202.1461 & \multicolumn{2}{c|}{}& \pur{170118.4611}\\
504 && 41283.0173& 41283.0173& 41283.0173& 41283.0173& 41283.0173 & \multicolumn{2}{c|}{}& \pur{172312.0002}\\\cline{3-3}
&\multicolumn{2}{c|}{}&..&..&..&..&\multicolumn{2}{c|}{}&..\\
761 & \multicolumn{2}{c|}{}& 64851.1221& 64851.1221& 64851.1221& 64851.1221 & \multicolumn{2}{c|}{}& \pur{1153652.4623}\\
762 & \multicolumn{2}{c|}{}& 64912.1209& 64912.12\underline{10}& 64912.1209& 64912.1209 & \multicolumn{2}{c|}{}& \pur{1162563.6568}\\
763 & \multicolumn{2}{c|}{}& 64983.7916& 64983.791\underline{8}& 64983.7916& 64983.7916 & \multicolumn{2}{c|}{}& \pur{1170839.8358}\\\cline{5-5}
&\multicolumn{2}{c|}{}&..&&..&..&\multicolumn{2}{c|}{}&..\\
955 & \multicolumn{2}{c|}{}& 82573.2725 && 82573.2725& 82573.2725 & \multicolumn{2}{c|}{}& \pur{3442594.0097}\\
956 & \multicolumn{2}{c|}{}& 82613.2799 && 82613.2799& 82613.2799 & \multicolumn{2}{c|}{}& \pur{3444776.2656}\\
957 & \multicolumn{2}{c|}{}& 82772.9680 && 82772.9680& 82772.9680 & \multicolumn{2}{c|}{}& \pur{3444776.4683}\\\cline{6-6}\cline{10-10}
&\multicolumn{2}{c|}{}&..&\multicolumn{2}{c|}{}&..&\multicolumn{3}{c|}{}\\
1014 & \multicolumn{2}{c|}{}& 88010.7113 &\multicolumn{2}{c|}{}& 88010.7113 & \multicolumn{3}{c|}{}\\
1015 & \multicolumn{2}{c|}{}& 88090.9183 &\multicolumn{2}{c|}{}& 88090.9183 & \multicolumn{3}{c|}{}\\
1016 & \multicolumn{2}{c|}{}& 88095.6490 &\multicolumn{2}{c|}{}& 88095.6490 & \multicolumn{3}{c|}{}\\
&\multicolumn{2}{c|}{}&..&\multicolumn{2}{c|}{}&..&\multicolumn{3}{c|}{}\\
\hline
\end{tabular}}
\end{table}

The results of the bridge behave, as expected, similar to the results of the frame structure and matching tendencies can be observed among the three different ansatz functions and the FEM results. In both cases, the NAT using the \textit{alternative ansatz} delivers outstanding results and outperforms the conventional ansatz functions and the FEM computations in many ways.
\section{Conclusion}
In the present paper, we extended the NAT to the computation of natural frequencies of plane frame structures based on the Euler-Bernoulli beam theory. Since analytic solutions are used, no discretization refinement is necessary and exact results are obtained.
We have formally introduced a frame structure as a set of nodes, beams, bearings, springs, and external loads. This allows us to formulate the boundary and interface conditions in a systematic way for an arbitrary system. In combination with the \emph{alternative ansatz}, the implemented method provides stable computations of natural frequencies also for higher frequency ranges and complicated structures with multiple beams. 

In future work, various extensions of the presented method should be easily possible. This includes the consideration of concentrated masses, dampers, and cracked beams in the two-dimensional setting. Furthermore, the extension to space frames and non-linear beam geometries would be an interesting research direction. Implementing the Timoshenko beam theory in the existing code can also be considered.



\appendix
\section{Appendix}

\begin{table}[H]
	\centering\setlength{\tabcolsep}{10pt}
	\makegapedcells
	\caption{A list of the different support types and its support conditions.}
	\label{tab:supports}
	\resizebox{7cm}{!}{
	\begin{tabular}{|c|c|c|} \hline
		\shortstack{name of \\ support type} 	&	\shortstack{symbolic \\ representation}	& \shortstack{supporting \\ conditions} \\ \hline 
		\rule{0pt}{0ex}	Pinned  & \begin{tikzpicture}
			\draw[ thick, ->,rotate around={0:(-0.5,0.5)}] (-0.5,0.5)--(0,0.5) node [right]{\footnotesize $\zeta$};
			\draw[ thick, ->,rotate around={-90:(-0.5,0.5)}] (-0.5,0.5)--(0,0.5) node [ below]{\footnotesize $\eta$};
			\draw (-0.5,0.5)circle(1pt);
			\coordinate (a) at (.5,0);
			\coordinate (b) at (1,0);
			\def\offset{.1};
			\draw [ thick](a)--(b);
			\draw [ thick, dashed]($(a)+(\offset,-\offset)$)--($(b)+(0,-\offset)$);
			\Schnittufer{b}[0][.5];	
			\begin{scope}[shift={(.5,0)},rotate=0,scale=1]   
				\draw[thick] (0,0) -- (-0.25,-0.4) -- (0.25,-0.4) -- (0,0);
				\draw[thick] (-0.36,-0.4) -- +(0.72,0.);
				\def\step{0.1}
				\def\stepp{0.4}
				\foreach \x in {0.,\step,...,.3}{
					\draw[thin] (-.35+\x,-.4) -- (-.35,-.4-\x);}
				\foreach \x in {.3,\stepp,...,.69}{
					\draw[thin] (-.35+\x,-.4) -- (-.65+\x,-.7);}
				\foreach \y in {0.,\step,...,.3}{
					\draw[thin] (.35,-.4-\y) -- (.05+\y,-.7);}   
			\end{scope}
		\end{tikzpicture}	& \shortstack{$ M = 0 $ \\ $ U_\zeta = 0 $ \\ $ U_\eta = 0 $}
		\\[1ex] \hline
		Roller & \begin{tikzpicture}
			\coordinate (a) at (0,0);
			\coordinate (b) at (.5,0);
			\def\offset{.1};						
			\draw [ thick](a)--(b);
			\draw [ thick, dashed]($(a)+(\offset,-\offset)$)--($(b)+(0,-\offset)$);
			\Schnittufer{b}[0][.5];
			\begin{scope}[shift={(a)},rotate=0,scale=1] 
				\draw[thick] (0,0) -- (-0.25,-0.4) -- (0.25,-0.4) -- (0,0);
				\draw[thick] (-0.36,-0.4) -- +(0.72,0.);
				\draw[thick] (-0.36,-0.5) -- +(0.72,0.);
				\def\step{0.1}
				\def\stepp{0.4}
				\foreach \x in {0.,\step,...,.3}{
					\draw[thin] (-.35+\x,-.5) -- (-.35,-.5-\x);}
				\foreach \x in {.3,\stepp,...,.69}{
					\draw[thin] (-.35+\x,-.5) -- (-.65+\x,-.8);}
				\foreach \y in {0.,\step,...,.3}{
					\draw[thin] (.35,-.5-\y) -- (.05+\y,-.8);}         
			\end{scope}	
		\end{tikzpicture}	& \shortstack{$ M = 0 $ \\ $ S_\zeta = 0 $ \\ $ U_\eta = 0 $}
		\\ \hline
		Clamped & \begin{tikzpicture}
			\coordinate (a) at (0,0);
			\coordinate (b) at (0,.5);
			\draw [ thick](a)--(b);
			\def\offset{.1};									
			\draw [ thick, dashed]($(a)+(\offset,0)$)--($(b)+(\offset,0)$);
			\Einspannung{a};
			\Schnittufer{b}[90][.5];				
		\end{tikzpicture}	& \shortstack{$ \psi = 0 $ \\ $ U_\zeta = 0 $ \\ $ U_\eta = 0 $}
		\\ \hline
		Parallel guide & \begin{tikzpicture}
			\coordinate (a) at (0,0);
			\coordinate (b) at (0,0.5);
			\def\offset{.1};											
			\draw [ thick](a)--(b);
			\draw [ thick, dashed]($(a)+(\offset,0)$)--($(b)+(\offset,0)$);
			\Parallelfuehrung{a};
			\Schnittufer{b}[90][.5];								
		\end{tikzpicture}	& \shortstack{$ \psi = 0 $\\$ S_\zeta = 0 $ \\ $ U_\eta = 0 $  }
		\\ \hline
	\end{tabular}}
\end{table}
\bibliography{mybibfile}
\end{document}